\documentclass{article}       
\usepackage{graphicx}
\usepackage{url} 
\usepackage{amssymb}
\usepackage{amsmath, amsthm}
\usepackage{mathptmx}      

\newtheorem{theorem}{Theorem}
\newtheorem{lemma}{Lemma}

\begin{document}

\title{A variation of a congruence of Subbarao for $n=2^{\alpha}5^{\beta}$\thanks{The author is supported by Croatian Science Foundation grant number 6422.}
}

\author{Sanda Buja\v{c}i\'{c}  \footnote{Department of Mathematics 
		Radmile Matej\v{c}i\'{c} 2, 51000 Rijeka, Croatia 
		Tel.: +385-51-584654
		Fax: +385-51-584699
		sbujacic@math.uniri.hr }
}

\date{}



\maketitle

\begin{abstract}
There are many open problems concerning the characterization of the positive integers $n$ fulfilling certain congruences and involving the 
Euler totient function $\varphi$ and the sum of positive divisors function $\sigma$ of the positive integer $n$. In this work, we deal with the congruence of the form

$$\label{cong}
n\varphi(n)\equiv2\pmod{\sigma(n)}
$$
and we prove that the only positive integers of the form $2^{\alpha}5^{\beta}, \enspace \alpha, \beta\geq0,$ that satisfy the above congruence are $n=1, 2, 5, 8$.\\
\textbf{Keywords}\\
Euler's totient function  sum of positive divisors  Pellian equations  congruence of Subbarao
\end{abstract}

\section{Introduction} \label{intro}

Wilson's theorem is a well known characterization of prime 
numbers. It states that a positive integer $p>1$ is a prime number if and only 
if the congruence of the form $$(p-1)!+1\equiv0\pmod{p}$$ is satisfied. There is probably no other so simple characterization of prime 
numbers in the form of a congruence, but there are many open problems concerning the characterization of the 
positive integers fulfilling certain congruences and involving functions 
$\varphi$ and $\sigma$, where $\varphi(n)$ and $\sigma(n)$ stand for the 
Euler totient function and the sum of positive divisors function of the 
positive integer $n$, respectively. 

In 1932 D. H. Lehmer \cite{lehmer} was dealing with the congruence of the form 
\begin{equation}\label{eq:lehmer}
n-1\equiv0\pmod{\varphi(n)}.
\end{equation}
This problem is known as Lehmer's totient problem. Despite the fact 
that the congruence $(\ref{eq:lehmer})$ is satisfied by every prime number, Lehmer's totient problem is an open problem because it is still not known whether there exists a composite number that satisfies it. Lehmer proved that, if there 
exists a composite number that satisfies the congruence (\ref{eq:lehmer}), then it must be odd, square free and it must have 
at least seven distinct prime factors. In 1944, F. Schuh \cite{schuh}
improved Lehmer's result and showed that such composite number must have at 
least eleven distinct prime factors. The best current result shows that, if such composite number exists, it has to have at least fourteen distinct prime factors \cite{cohen}.

M. V. Subbarao was considering the congruence of the form 
\begin{equation}\label{eq:subbarao}
n\sigma(n)\equiv 2\pmod{\varphi(n)}.
\end{equation}
He proved \cite{subbarao} that the only composite numbers that 
satisfy the congruence (\ref{eq:subbarao}) are numbers $4, 6$ and $22$. 

A. Dujella and F. Luca were dealing with the congruence of the form
\begin{equation}\label{eq:version}
n\varphi(n)\equiv2\pmod{\sigma(n)},
\end{equation}
which is a variation of the congruence (\ref{eq:subbarao}). They have proved
\cite{luca} that there are only finitely many positive integers that 
satisfy the congruence (\ref{eq:version}) and whose prime factors belong to a fixed finite set. They have proved that when this finite set consists only of two primes $2$ and $3$, then the only positive integers of the form $n=2^a3^b, \ a,  b \geq 0,$ that satisfy the congruence (\ref{eq:version}) are $n=1, 2, 3, 8, 9.$

We deal with the variation of the congruence of Subbarao for $n=2^{\alpha}5^{\beta}, \alpha, \beta\geq0$. The main result of this paper is the following theorem.

\begin{theorem}\label{main_theo}
	The only positive integers $n$ of the form $n=2^{\alpha}5^{\beta}, \ \alpha, \beta \geq 0, $ that satisfy the congruence
	$$n\varphi(n) \equiv 2 \pmod{\sigma(n)}$$ are  numbers $n=1, 2, 5, 8.$
\end{theorem}

\section{The proof of Theorem \ref{main_theo}}

\subsection{The cases $\alpha=0$ and $\beta=0$}

	It is easily seen that all the prime numbers satisfy the congruence 
	(\ref{eq:version}). Let $p$ be a prime. In this case, we have $\varphi(p)=p-1$ 
	and $\sigma(p)=p+1$. The variation of the congruence of Subbarao (\ref{eq:version}) 
	becomes
	$$p(p-1)-2=(p+1)(p-2)\equiv0\pmod{(p+1)}.$$
	The congruence (\ref{eq:version}) is 
	satisfied for all the prime numbers, or more precisely, 
	$$p(p-1)\equiv2\pmod{(p+1)}.$$ Hence, the prime numbers $2$ and $5$ satisfy the congruence (\ref{eq:version}). 
	
	The remaining part of the proof deals with the composite numbers of the form 
	$n=2^{\alpha}5^{\beta}, \ \alpha, \beta\geq0$. 
	
	For start, let $\beta=0$ which implies dealing with 
	the positive integers of the form $n=2^{\alpha}, \enspace \alpha\geq2$. We define 
	$$D:=\sigma(2^{\alpha})=2^{\alpha+1}-1.$$ We may notice $2^{\alpha+1}\equiv1\pmod{D}$. 
	Because of $(\ref{eq:version})$, we obtain 
	$$2^{\alpha}\cdot2^{\alpha}\left(1-\frac{1}{2}\right)\equiv2\pmod{D},$$
	$$2^{2(\alpha+1)}\equiv2^4\pmod{D},$$
	$$(2^{\alpha+1}-1)(2^{\alpha+1}+1)-15\equiv0\pmod{D}.$$
	The condition $D\mid ((2^{\alpha+1}-1)(2^{\alpha+1}+1)-15)$ is satisfied if and only if 
	$D\mid 15$, or more precisely, if and only if $$(2^{\alpha+1}-1)\mid 15.$$ For 
	$\alpha\geq2$, the condition $(2^{\alpha+1}-1)\mid 15$ is satisfied only when $\alpha=3$. So, 
	$n=2^3$ is the only positive integer of the form $n=2^\alpha, \ \alpha \geq2,$ that  
	satisfies the variation of the congruence of Subbarao $(\ref{eq:version})$.
	
	Now, let $\alpha=0$. We deal with the positive integers of the form $n=5^{\beta}, \enspace 
	\beta\geq2$. We define $$D:=\sigma(5^{\beta})=\frac{5^{\beta+1}-1}{4}.$$ As in the previous 
	case, it is easy to notice that $5^{\beta+1}\equiv1\pmod{D}$. Because of $(\ref{eq:version})$, we obtain 
	$$5^{2\beta-1}\cdot2^2\equiv2\pmod{D},$$ $$5^{2(\beta+1)}\cdot2^2\equiv5^3\cdot 
	2\pmod{D}.$$
	Using $5^{\beta+1}\equiv1\pmod{D}$, the previous congruence implies  
	$D\mid 246$, which is not possible for $\beta\geq2$. Consequently, there are no positive integers of 
	the form $n=5^{\beta}, \enspace \beta\geq2,$ that satisfy the congruence 
	$(\ref{eq:version})$.
	
	\subsection{Useful congruences}
	
	The remaining part of the proof deals with the most general case, or more 
	precisely, with the positive integers of the form $n=2^{\alpha}5^{\beta}, \ \alpha, \beta\in \mathbb{N}.$ We start by defining $M:=2^{\alpha+1}-1$ and $N:=\frac{5^{\beta+1}-1}{4}$. As in the 
	previous cases, we use congruences $$2^{\alpha+1}\equiv1\pmod{M}$$ and 
	$$5^{\beta+1}\equiv1\pmod{N}.$$
	Supposing the congruence $(\ref{eq:version})$ is satisfied, we get
	\begin{equation}\label{eq:t}
	2^{2\alpha+1}\cdot5^{2\beta-1}\equiv2\pmod{MN}.
	\end{equation}
	\noindent Multiplying (\ref{eq:t}) by $2\cdot5^3,$ we obtain
	$$2^{2(\alpha+1)}\cdot5^{2(\beta+1)}\equiv500\pmod{MN}.$$
	\indent Since $2^{\alpha+1}\equiv1\pmod{M},$ we get  
	$5^{2(\beta+1)}\equiv500\pmod{M}$. Analogously, because of $5^{\beta+1}\equiv1\pmod{N}$, we conclude $2^{2(\alpha+1)}\equiv500\pmod{N}$. 
	
	For $M\mid (2^{\alpha+1}-1)$, we have $M\mid (2^{2(\alpha+1)}-1)$. Similarly, $N\mid (5^{2(\beta+1)}-1)$. Hence, we get \begin{equation}\label{eq:21}
	M, N \mid (2^{2(\alpha+1)}+5^{2(\beta+1)}-501).
	\end{equation}
	
	Our next step is to show that $\alpha$ and $\beta$ are even and $M$ and $N$ are coprime. Let $G:=\textnormal{gcd}(M, N)$. Because $2^{\alpha+1}\equiv1\pmod{M}$ and $5^{\beta+1}\equiv1\pmod{N},$ we get $$2^{\alpha+1}\equiv 
	5^{\beta+1}\equiv1\pmod{G}.$$ Because of $(\ref{eq:21})$, we conclude $G\mid -499$. 
	Number $499$ is a prime, so $G=1$ or $G=499$.  For start, we can assume 
	$G=499$. In this case, we know that $499\mid M$, or, more precisely, $$499\mid 
	\left(2^{\alpha+1}-1\right).$$ The order of $2$ modulo $499$ is $166$, so $166\mid 
	(\alpha+1)$. Especially, $2\mid (\alpha+1)$. Hence, $\alpha$ is odd. We can notice 
	that $M$ can be expressed as $M=2^{\alpha+1}-1=2^{2k}-1,$ for $k\in\mathbb{N}$. 
	Obviously, $3\mid M$. Hence, $3\mid (n\varphi(n)-2)$, or, more precisely, $3\mid 
	(2^{2\alpha+1}\cdot5^{2\beta-1}-2)$, which is not possible. As a consequence, we conclude $499\nmid M$ which implies that $G=1$. We have proved that $M, N$ are coprime. 
	
	The next step is to determine the parity of $\alpha$ and $\beta$. For start, we assume that $\alpha$ is odd which implies that $\alpha+1$ is even. So,
	$$M=2^{\alpha+1}-1=2^{2k}-1,$$ for $k\in\mathbb{N}$. Obviously, $3\mid M$ and according to our hypothesis, $3\mid 
	(2^{2\alpha+1}\cdot5^{2\beta-1}-2)$, which is not possible. We conclude that $\alpha$ is even.
	
	Now, we assume that $\beta$ is odd. In that case, we write $$5^{\beta+1}-1=5^{2k}-1,$$ for 
	$k\in\mathbb{N}$. Obviously, $24\mid (5^{2k}-1)$ and because $6\mid N$ and 
	$N\mid (2^{2\alpha+1}\cdot5^{2\beta-1}-2)$, we get $6\mid 
	(2^{2\alpha+1}\cdot5^{2\beta-1}-2)$, which is not possible. Hence, $\beta$ is even, too.
	
	Hence, we have proved 
	that $M$ and $N$ are odd and coprime numbers.
	
	\indent As a consequence of (\ref{eq:21}), we may notice $$MN\mid 
	(2^{2(\alpha+1)}+5^{2(\beta+1)}-501).$$ On the other hand, 
	$$4MN=(2^{\alpha+1}-1)(5^{\beta+1}-1),$$
	and obviously $2^{2(\alpha+1)}+5^{2(\beta+1)}-501\equiv0\pmod{4}$.\\
	\indent Let $x:=2^{\alpha+1}$ and $y:=5^{\beta+1}$. The initial 
	problem is now represented by the equation of the form
	\begin{equation}\label{eq:22}
	x^2+y^2-501=c(x-1)(y-1),
	\end{equation} for some $c\in\mathbb{N}$. 
	
	Since numbers $\alpha$ and $\beta$ are even, the following congruences hold $$x\equiv0\pmod{8}, \enspace x^2\equiv0\pmod{8},$$
	$$y\equiv5\pmod{8}, \enspace y^2\equiv1\pmod{8}.$$ From $(\ref{eq:22})$, we get $4c\equiv4\pmod{8}$, which is satisfied for
	\begin{equation}\label{eq:kongr1}
	c\equiv1\pmod{2}.
	\end{equation} 
	
	We also notice that congruences $$x\equiv2\pmod{3},\enspace x^2\equiv1\pmod{3},$$
	$$y\equiv2\pmod{3},\enspace y^2\equiv1\pmod{3}$$ are satisfied. From 
	$(\ref{eq:22})$, we easily get
	\begin{equation}\label{eq:kongr2}
	c\equiv2\pmod{3}.
	\end{equation}

	We conclude $$x\equiv3\pmod{5},\enspace x^2\equiv4\pmod{5} \enspace \enspace
	\textnormal{for} \enspace \alpha\equiv2\pmod{4}, $$
	$$x\equiv2\pmod{5}, \enspace x^2\equiv4\pmod{5} \enspace \enspace \textnormal{for} 
	\enspace \alpha\equiv0\pmod{4}.$$
	Obviously, $$y\equiv y^2\equiv0\pmod{5}.$$ Bringing everything together, we obtain that
	$$c\equiv1\pmod{5},\enspace \textnormal{for} \enspace \alpha\equiv2\pmod{4},$$
	or
	\begin{equation}\label{eq:kongr3}
	c\equiv2\pmod{5},\enspace \textnormal{for} \enspace \alpha\equiv0\pmod{4}.
	\end{equation}

	Now, we try to determine which of the above residue classes modulo $5$ is satisfied by the number $c$ that is introduced in our problem. 
	
	Let $t=2^{\alpha}\cdot5^{\beta-1}$. We get that
	$$5t^2=2^{2\alpha}\cdot5^{2\beta-1}.$$ According to $(\ref{eq:t}),$ we conclude
	$5t^2\equiv1\pmod{M}$, which implies  
	$\left(\frac{5}{M}\right)=\left(\frac{M}{5}\right)=1$. In this case
	$M\equiv1, 4\pmod{5}$. Since $M=2^{\alpha+1}-1$, we get $2^{\alpha+1}-1\equiv 
	1\pmod{5}$ or $2^{\alpha+1}-1\equiv 4\pmod{5}$. The first congruence is satisfied for $\alpha\equiv0\pmod{4},$ while the second congruence, $2^{\alpha+1}-1\equiv 4\pmod{5}$, is never satisfied. Consequently, we consider only positive integers $c$ that satisfy the congruence $$c\equiv2\pmod{5}.$$
	
	Taking into account congruences $(\ref{eq:kongr1}), (\ref{eq:kongr2})$ and 
	$(\ref{eq:kongr3})$ and using Chinese Remainder Theorem, we determine that required positive integers $c$ satisfy
	\begin{equation}\label{eq:kongr}
	c\equiv17\pmod{30}.
	\end{equation}
	
	\subsection{Pellian equations}
	\indent We "diagonalize" the equation $(\ref{eq:22})$ as in 
	\cite{luca}. Let \begin{equation}\label{eq:23}
	X:=cy-c-2x,
	\end{equation}
	\vspace{-20pt}
	\begin{equation}\label{eq:24}
	Y:=cy-c-2y.
	\end{equation}
	Then 
	$$(c+2)Y^2-(c-2)X^2-(-1996c+4008)=-4(c-2)(x^2+y^2-501-c(x-1)(y-1))=0.$$
	This method has resulted with the Pellian equation of the form
	\begin{equation}\label{eq:25}
	(c+2)Y^2-(c-2)X^2=-1996c+4008.
	\end{equation}

	\indent Let $X=0$. In this case, the Pellian equation
	(\ref{eq:25})
	becomes $$Y^2=\frac{-1996c+4008}{c+2}.$$ The only integer solution of the above equation is $Y=\pm2$ for $c=2$. Since $c=2$ does not satisfy the congruence (\ref{eq:kongr}), in our case $Y=\pm2$ is not the solution of the equation (\ref{eq:25}).\\
	\indent Let $Y=0$. The initial Pellian equation (\ref{eq:25}) becomes
	$$X^2=\frac{1996c-4008}{c-2}.$$
	The right-hand side of the equation is an integer for $c=1, 3, 4, 6, 
	10, 18$. Those numbers do not satisfy the congruence (\ref{eq:kongr}). 
	The right-hand side of the above equation is not a perfect square for such positive integers $c$, so we conclude there does not exist a solution $X$ of the Pellian equation (\ref{eq:25}).
	
	\indent Now we deal with the general case. Let $(X, Y)$ be a solution of the 
	equation (\ref{eq:25}) in positive integers. In this case, $\frac{X}{Y}$ is a 
	good rational approximation of the irrational number $\sqrt{\frac{c+2}{c-2}}$. More 
	precisely,
	$$\left| \frac{X}{Y}-\sqrt{\frac{c+2}{c-2}}\right| 
	=\frac{1996c-4008}{(\sqrt{c+2}Y+\sqrt{c-2}X)\sqrt{c-2}Y}\leq\frac{1996(c-2)}{
		\sqrt{c^2-4}Y^2}<\frac{1996}{Y^2}.$$
	The rational approximation of the form
	\begin{equation}\label{eq:26}
	\left| \frac{X}{Y}-\sqrt{\frac{c+2}{c-2}}\right| <\frac{1996}{Y^2}
	\end{equation}
	is not good enough to conclude that $\frac{X}{Y}$ is a convergent of continued fraction expansion of
	$\sqrt{\frac{c+2}{c-2}}$. We use Worley and Dujella's theorem from \cite{worley} and 
	\cite{dujellarsa}.
	
	\begin{theorem}[Worley, Dujella]
		\textit{Let $\alpha$ be an irrational number and let $a, b\neq0$ be coprime nonzero integers satisfying the inequality $$\left| 
			\alpha-\frac{a}{b}\right| <\frac{H}{b^2},$$
			where $H$ is a positive real number. Then $$(a, b)=(rp_{m+1}\pm sp_m,\enspace 
			rq_{m+1}\pm sq_m),$$ for $m, r, s\in\mathbb{N}_0$ such that
			$rs<2H,$ where $\frac{p_m}{q_m}$ is $m$--th convergent from continued fraction expansion of irrational number $\alpha$.}
	\end{theorem}
	
	According to Worley and Dujella's theorem, we get that every solution $(X, Y)$ of the Pellian
	equation (\ref{eq:25}) is of the form
	$$X=\pm d(rp_{k+1}+ up_k), \enspace Y=\pm d(rq_{k+1}+ uq_k)$$
	for some $k\geq-1$, $u\in\mathbb{Z}$, $r$ nonnegative positive 
	integer and $d=\textrm{gcd}(X, 
	Y)$ for which the inequality $$|ru|<2\cdot\frac{1996}{d^2}$$ holds. 
	
	In order to determine all the integer solutions of the Pellian equation (\ref{eq:25}), we also use a lemma from 
	\cite{jadri}.
	\begin{lemma}[Dujella, Jadrijevi\'{c}]\label{jadrilema}
		Let $\alpha\beta$ be a positive integer which is not a perfect square and let
		$p_k/q_k$ be the $k$-th convergent of continued fraction expansion of
		$\sqrt{\frac{\alpha}{\beta}}$. Let the sequences $(s_k)_{k\geq-1}$ and $(t_k)_{k\geq-1}$ be
		the sequences of the integers appearing in the continued fraction expansion of $\frac{\sqrt{\alpha\beta}}{\beta}$. Then
		\begin{equation}\label{worley}
		\alpha(rq_{k+1}+uq_k)^2-\beta(rp_{k+1}+up_k)^2=(-1)^k(u^2t_{k+1}+2rus_{k+2}
		-r^2t_{k+2}).
		\end{equation}
	\end{lemma}
	
	Applying Lemma \ref{jadrilema}, it is easy to conclude that we obtain 
	\begin{equation}\label{jadrilema1}
	(c+2)Y^2-(c-2)X^2=d^2(-1)^k(u^2t_{k+1}+ 2rus_{k+2}-r^2t_{k+2}),
	\end{equation} where ${(s_k)}_{k\geq-1}$ and ${(t_k)}_{k\geq-1}$ are sequences 
	of integers appearing in the continued fraction expansion of the quadratic irrationality 
	$\sqrt{\frac{c+2}{c-2}}$. Our next step is to determine  the continued fraction expansion of
	$\sqrt{\frac{c+2}{c-2}},$ where $c$ is an odd positive integer.
	
	From the continued fraction expansion algorithm we get
	$$s_0=0, \enspace t_0=c-2, \enspace a_0=1,$$
	$$s_1=c-2, \enspace t_1=4, a_1=\frac{c-3}{2},$$
	$$s_2=c-4, \enspace t_2=2c-5, \enspace a_2=1,$$
	$$s_3=c-1, \enspace t_3=1, \enspace a_3=2c-2,$$
	$$s_4=c-1, \enspace t_4=2c-5, \enspace a_4=1,$$
	$$s_5=c-4, \enspace t_5=4, \enspace a_5=\frac{c-3}{2},$$
	$$s_6=c-2, \enspace t_6=c-2, \enspace a_6=2,$$
	hence $$\sqrt{\frac{c+2}{c-2}}=\Big[ 1;\overline{\frac{c-3}{2}, 1, 2c-2, 1, \frac{c-3}{2}, 2}\Big].$$
	The length $l$ of the period of the continued fraction expansion of $\sqrt{\frac{c+2}{c-2}}$ is $l=6$, so we consider the equation $(\ref{jadrilema1})$ for $k=0, 1, 2, 3, 4, 5$ and determine all the positive integers $c$ that satisfy the congruence $(\ref{eq:kongr})$. From (\ref{eq:25}) and (\ref{jadrilema1}) we get 
	\begin{equation}\label{eq:dod}
	d^2(-1)^k(u^2t_{k+1}+2rus_{k+2}-r^2t_{k+2})=-1996c+4008.
	\end{equation}
	Obviously, $d$ can be $d=1$ or $d=2,$ for all $k=0, 1, 2, 3, 4, 5$.
	
	Let $k=0.$ From $(\ref{eq:dod}),$ we obtain
	$$d^2(u^2t_1+2rus_2-r^2t_2)=-1996c+4008,$$
	\begin{equation}\label{eq:opis}
	d^2(4u^2+2(c-4)ru-r^2(2c-5))=-1996c+4008.
	\end{equation}
	First, we deal with the cases when $d=1$ and $d=2$ and check if it is possible that both sides in the above equation are identical for each such $d$.

	For $d=1$ we get the system of two equations
	\[ 
	\begin{cases}
	4u^2-8ru+5r^2=4008,\\
	2ru-2r^2=-1996,
	\end{cases}
	\]
	that does not have integer solutions. Analogously, for $d=2$ we get
	the system
	\[
	\begin{cases}
	4u^2-8ru+5r^2=1002,\\
	2ru-2r^2=-499,
	\end{cases}
	\]
	which also does not have any integer solution.
	
	Generally, for $k=0$ and for all values of $d$, we obtain from (\ref{eq:opis}) that the positive integer $c$ is of the form  
	\begin{equation}\label{eq:prvic}
	c=\frac{4008-4d^2u^2+8d^2ru-5d^2r^2}{1996+2d^2ru-2d^2r^2}.
	\end{equation}
	Our goal is to determine all positive integers $c$ that satisfy the congruence (\ref{eq:kongr}), that are of the form (\ref{eq:prvic}) and for which the triples $(d, r, u)$ satisfy the conditions $d\in\mathbb{N}, \enspace r\in\mathbb{N}, \enspace u\in\mathbb{Z}, \enspace u\neq0$ and the inequality 
	\begin{equation}\label{inequality}
	d^2|ru|<3992.
	\end{equation}
	It is useful to mention that the latter condition implies $d\leq63$.\\
	\indent  Before dealing with the general case, we analyze the case when $ru=0$ from which we obtain
	$$c=\frac{4008-4d^2u^2}{1996},$$ for $(r, u)=(0, u)$ and 
	$$c=\frac{4008-5d^2r^2}{1996-2d^2r^2},$$ for $(r, u)=(r, 0)$.

	The equations do not hold for a positive integer $c$, except for $c=2$. Such $c$ does not satisfy the congruence $(\ref{eq:kongr})$ which allows us to conclude that in our case there are no integer solutions of the above equations that derive from these special cases.\\
	\indent An algorithm for generating triples $(d, r, u)$ that satisfy the inequality (\ref{inequality}) is created. This algorithm plugs these triples $(d, r, u)$ into $(\ref{eq:prvic})$ and checks if positive integers $c$ satisfy the congruence $(\ref{eq:kongr})$. 
		
	For $k=0$ we get $$c\in\{17,227,497,647,857,2537,3107,4937\}.$$ For each such positive integer $c$ we obtain a Pellian equation of the form $(\ref{eq:25})$.
	
	For $k=1$ the equation $(\ref{eq:dod})$ becomes
	$$-d^2(u^2(2c-5)+2ru(c-1)-r^2)=-1996c+4008.$$
	\indent If we consider the case when both sides in the above equation are identical, for $d=1$ we get $$5u^2+2ru+r^2=4008, \enspace 2u^2+2ru=1996,$$ while for $d=2,$ we obtain  $$5u^2+2ru+r^2=1002, \enspace 2u^2+2ru=499.$$ There are no integer solutions of both systems.\\
	\indent Generally, from $(\ref{eq:dod})$ we get that $c$ is represented by
	$$c=\frac{5d^2u^2+2d^2ru+d^2r^2-4008}{2d^2u^2+2d^2ru-1996}.$$
	The described algorithm is used to get the following values for $c$: 
	$$c\in\{17, 227, 497, 647, 857, 2537, 3107, 4937\}.$$
	
	For $k=2$ we have
	$$d^2(u^2+2ru(c-1)-r^2(2c-5))=-1996c+4008.$$
	If both sides in the above equation are identical, for $d=1$ we obtain  $$u^2-2ru+5r^2=4008, \enspace 2ru-2r^2=-1996$$ which is the system of two equations that does not have integer solutions. Analogously, for $d=2$ the system 
	$$u^2-2ru+5r^2=1002, \enspace 2ru-2r^2=-499$$ has no integer solutions.\\
	Generally, the positive integer $c$ obtained from $(\ref{eq:dod})$ when $k=2$ is of the form
	$$c=\frac{d^2u^2-2d^2ru+5d^2r^2-4008}{2d^2r^2-2d^2ru-1996}.$$
	We get
	$$c\in\{17,227,497,647,857,2537,3107,4937\}.$$ 
		
	When $k=3$ from (\ref{eq:dod}) we get
	$$-d^2(u^2(2c-5)+2ru(c-4)-4r^2)=-1996c+4008.$$
	For $d=1$ and $d=2$ the following systems are obtained respectively
	$$5u^2+8ru+4r^2=4008, \enspace 2u^2+2ru=1996$$ and $$5u^2+8ru+4r^2=1002, \enspace 2u^2+2ru=499.$$ Like in previous cases, these systems do not have integer solutions.
	Generally, the positive integer $c$ is of the form
	$$c=\frac{5d^2u^2+8d^2ru+4d^2r^2-4008}{2d^2u^2+2d^2ru-1996}.$$
	We get the following values for $c$ in this case:
	$$c\in\{17,227,497,647,857,2537,3107,4937\}.$$
	
	Analogously, for $k=4$ we get
	$$d^2(4u^2+2ru(c-2)-r^2(c-2))=-1996c+4008.$$
	For $d=1$ we obtain $$4u^2-4ru+2r^2=4008, \enspace 2ru-r^2=1996,$$ while for $d=2$ we get $$4u^2-4ru+2r^2=1002, \enspace 2ru-r^2=499.$$ Both systems do not have integer solutions. Generally,
	$$c=\frac{4d^2u^2-4d^2ru+2d^2r^2-4008}{d^2r^2-2d^2ru-1996}.$$
	For $k=4$ we get $$c\in\{17,227,497,647,857,2537,3107,4937\}.$$
	
	Finally, for $k=5$ from (\ref{eq:dod}) we get
	$$-d^2(u^2(c-2)-r^2(c-2))=-1996c+4008.$$
	If we take into account the case when both sides of the above equation are identical, for $d=1$ we obtain the following system of equations $$2u^2-2r^2=4008, \enspace r^2-u^2=-1996,$$ and for $d=2$ we get $$2u^2-2r^2=1002, \enspace r^2-u^2=499.$$ Both systems do not have integer solutions. Generally,
	\begin{equation}\label{eq:c}c=\frac{2d^2u^2-2d^2r^2-4008}{d^2u^2-d^2r^2-1996}.\end{equation}
	There is no $c$ of the form (\ref{eq:c}) which satisfies the congruence $(\ref{eq:kongr})$.
	
	We gather all the possible positive integers $c\equiv17 \pmod{30}$ for $k=0, 1, 2, 3, 4, 5$ that we have determined using described algorithm and set a Pellian equation of the form $(\ref{eq:25})$ for every obtained $c$. The Pellian equations are
	$$19Y^2-15X^2=-29924, \enspace \textrm{for} \enspace c=17,$$
	$$229Y^2-225X^2=-449084, \enspace \textrm{for} \enspace c=227,$$
	$$499Y^2-495X^2=-988004, \enspace \textrm{for} \enspace c=497,$$
	$$649Y^2-645X^2=-1287404, \enspace \textrm{for} \enspace c=647,$$
	$$859Y^2-855X^2=-1706564, \enspace \textrm{for} \enspace c=857,$$
	$$2539Y^2-2535X^2=-5059844, \enspace \textrm{for} \enspace c=2537,$$
	$$3109Y^2-3105X^2=-6197564, \enspace \textrm{for} \enspace c=3107,$$
	$$4937Y^2-4935X^2=-9850244, \enspace \textrm{for} \enspace c=4937.$$
	In order to determine whether these Pellian equations have solutions, we use \cite{alpern}.
	
	\indent First of all, we assume that $X, Y$ are of the form $(\ref{eq:23})$, $(\ref{eq:24})$, respectively. We can easily determine that $X$ satisfies the following congruences
	$$X\equiv0\pmod{4}, \enspace X\equiv1\pmod{3}, \enspace X\equiv4\pmod{5}, \enspace \textrm{hence} \enspace X\equiv4\pmod{60}.$$ We set $X=60i+4, \enspace i\in\mathbb{Z}$.
	
	Analogously, $$Y\equiv2\pmod{4}, \enspace Y\equiv1\pmod{3}, \enspace Y\equiv3\pmod{5},$$ hence, $Y\equiv58\pmod{60}$. We set $Y=60j+58, \enspace j\in\mathbb{Z}$. 
	
	We deal with the first Pellian equation $$19Y^2-15X^2=-29924.$$ For $X=60i+4$ the above equation becomes
	$$19Y^2-54000i^2-7200i+29684=0.$$
	Using \cite{alpern} we determine that this Pellian equation does not have any integer solution.

	The next Pellian equation is $$229Y^2-225X^2=-449084.$$ 
	For $Y=60j+58,\enspace j\in\mathbb{Z}$, this equation becomes
	$$824400j^2-225X^2+1593840j+1219440=0,$$
	and it does not have any integer solution according to \cite{alpern}. 
	
	The next Pellian equation is \begin{equation}\label{eq:treca}499Y^2-495X^2=-988004.\end{equation}
	
	We can notice that the equation (\ref{eq:treca}) has integer solutions for $X\equiv4\pmod{60}$ and $Y\equiv58\pmod{60}$. We need to get some additional conditions for $X, Y$ in order to reach the conclusion that the equation (\ref{eq:treca}) does not have any integer solution for such $X, Y$.\\
	\indent Additionally, we know that $$Y=cy-c-2y=c(y-1)-2y\equiv-2\pmod{(c-2)}.$$ In this case, we have $c-2=495=3^2\cdot5\cdot11$, which implies $Y\equiv-2\pmod{3^2\cdot5\cdot11}$, or, more precisely, $$Y\equiv-2\pmod{9}, \enspace Y\equiv-2\pmod{5}, \enspace Y\equiv-2\pmod{11}.$$ We already know $Y\equiv2\pmod{4}$, so we can easily get $$Y=c(y-1)-2y=497(y-1)-2y\equiv-2y\equiv21, 28, 34, 61, 69\pmod{71}.$$ We set one Pellian equation of the form (\ref{eq:25}) for each residue that we get after dividing $Y$ by $71$ and we analyze each of these equations.
	
	We have $$Y\equiv2\pmod{4}, \enspace Y\equiv 7\pmod{3^2}, \enspace Y\equiv9\pmod{11}, \enspace Y\equiv21\pmod{71}.$$ We get $Y\equiv11878\pmod{140580}$, hence $Y=140580j+11878, \enspace j\in\mathbb{Z}$. For such $Y$ the equation (\ref{eq:treca}) becomes $$9861605463600 j^2 - 495 X^2 + 1666469621520 j + 70403343120=0.$$
	According to \cite{alpern} it does not have any integer solution.\\
	\indent For $$Y\equiv2\pmod{4}, \enspace Y\equiv 7\pmod{3^2}, \enspace Y\equiv9\pmod{11}, \enspace Y\equiv28\pmod{71},$$ we conclude $Y\equiv27718\pmod{140580}$ and $Y=140580j+27718, \enspace j\in\mathbb{Z},$ so the equation (\ref{eq:treca}) becomes 
	$$9861605463600 j^2- 495 X^2 + 3888803247120 j+383376462480=0.$$
	Using \cite{alpern} we conclude that the above equation does not have any integer solution.
	
	In the case when
	$$Y\equiv2\pmod{4}, \enspace Y\equiv 7\pmod{3^2}, \enspace Y\equiv9\pmod{11}, \enspace Y\equiv34\pmod{71},$$ we get that $Y\equiv61387\pmod{140580}$. For $Y=140580j+61387, \enspace j\in\mathbb{Z},$ we get
	$$9861605463600j^2-495X^2+8612524891080 j+1880414508735 =0.$$
	This equation does not have any integer solution according to \cite{alpern}.
	
	For $$Y\equiv2\pmod{4}, \enspace Y\equiv 7\pmod{3^2}, \enspace Y\equiv9\pmod{11}, \enspace Y\equiv61\pmod{71},$$ we obtain $Y\equiv1978\pmod{140580}$, hence $Y=140580j+1978, \enspace j\in\mathbb{Z}$. We get the Pellian equation
	$$9861605463600 j^2- 495 X^2 + 277511105520 j+1953317520 =0.$$
	Using online calculator \cite{alpern} we determine that the above equation also does not have any integer solution.\\
	\indent Finally, for $$Y\equiv2\pmod{4}, \enspace Y\equiv 7\pmod{3^2}, \enspace Y\equiv9\pmod{11}, \enspace Y\equiv69\pmod{71},$$ we get $Y\equiv140578\pmod{140580}$, which we can write as $Y=140580j+140578, j\in\mathbb{Z}.$ For such $Y$ we get the Pellian equation 
	$$9861605463600 j^2 - 495 X^2+ 19722930329520 j+ 9861325855920=0$$
	which does not have any integer solution according to \cite{alpern}.\\
	\indent We have proved that, in our case, the Pellian equation (\ref{eq:treca}) does not have any integer solution.\\
	\indent The next Pellian equation is $$649Y^2-645X^2=-1287404.$$
	
	For $Y=60j+58,\ j\in\mathbb{Z},$ we get
	$$2336400j^2-645X^2+4517040j+3470640=0,$$
	which is the Pellian equation that does not have any integer solution according to \cite{alpern}.\\
	\indent The next Pellian equation is  $$859Y^2-855X^2=-1706564.$$
	
	For $X=60i+4,\ i\in\mathbb{Z},$ we get
	$$859Y^2-3078000i^2-410400i+1692884=0.$$
	By \cite{alpern}, the Pellian equation does not have any integer solution.\\
	\indent The next Pellian equation is $$2539Y^2-2535X^2=-5059844.$$ 
	It is known that $Y\equiv-2\pmod{(c-2)}$. In our case, we have $Y\equiv-2\pmod{2535}$, or precisely, $Y\equiv-2\pmod{3\cdot5\cdot11^2}$. We can conclude 
	$$Y\equiv1\pmod{3}, \enspace Y\equiv3\pmod{5}, \enspace Y\equiv167\pmod{169}.$$ It is already known from before that $Y\equiv2\pmod{4}$. We get that $Y\equiv10138\pmod{10140}$, or $Y=10140j+10138,\ j\in\mathbb{Z}$. The Pellian equation of the form
	$$261058964400 j^2-2535 X^2 + 522014946960j +  260961052560 =0, $$
	does not have any integer solution according to \cite{alpern}.\\
	\indent The penultimate Pellian equation is $$3109Y^2-3105X^2=-6197564.$$
	For $X=60i+4,\ i\in\mathbb{Z},$ we obtain $$3109Y^2-11178000i^2-1490400i+6147884=0.$$
	By \cite{alpern}, this Pellian equation does not have any integer solution.\\
	\indent The last Pellian equation
	$$4937Y^2-4935X^2=-9850244$$
	does not have any integer solution according to \cite{alpern}.
	
	Since Pellian equations of the form (\ref{eq:25}) obtained for all the possible values of positive integers $c$ that satisfy the congruence (\ref{eq:kongr}) do not have solutions $X, Y$ in positive integers, we conclude that there do not exist positive integers of the form $n=2^{\alpha}5^{\beta}, \ \alpha, \beta \in\mathbb{N},$ that satisfy the variation of the congruence of Subbarao (\ref{eq:version}). Consequently, the only positive integers of the form  $n=2^{\alpha}5^{\beta}, \ \alpha, \beta\geq0,$ that satisfy the congruence (\ref{eq:version}) are $n=1, 2, 5, 8.$\\

\noindent\textbf{\large{Acknowledgements}}\\
We would like to thank Professor Andrej Dujella for many valuable suggestions and his help with the preparation of this article and to Professor Andrzej Schinzel for valuable remarks.


\begin{thebibliography}{}
%
%
\bibitem{alpern} D.~Alpern, Quadratic Diophantine equation solver, \url{http://www.alpertron.com.ar}
\bibitem{cohen} G. L. Cohen, and P. Hagis Jr., On the Number of Prime Factors of $n$ is $\phi(n)\mid(n-1)$, Nieuw Arch. Wisk. \textbf{28},  177--185, (1980).
\bibitem{dujellarsa} A.~Dujella, Continued fractions and RSA with small secret exponents, Tatra Mt. Math. Publ. \textbf{29}, 101--112 (2004).
\bibitem{jadri} A.~Dujella and B.~Jadrijevi\'{c}, A family of quartic Thue inequalities,  Acta Arith. \textbf{111}, 61--76, (2004).
\bibitem{luca} A.~Dujella and F.~Luca, On a variation of a congruence of Subbarao,  J. Aust. Math. Soc. \textbf{93}, 85--90, (2012).
\bibitem{lehmer} D.~H.~Lehmer, On Euler's totient function, Bull. Amer. Math. Soc. \textbf{38}, 745--751, (1932).
\bibitem{schuh} F.~Schuh, Do there exist composite numbers $m$ for which $\varphi(m)\mid(m-1)$ (Dutch), Mathematica Zupten B,  \textbf{13}, 102--107, (1944).
\bibitem{subbarao} M.~V.~Subbarao, On two congruences for primality, Pacific J. Math.  \textbf{54}, 261--268, (1974).
\bibitem{worley} R.~T.~Worley, Estimating $|\alpha-\frac{p}{q}|$, J. Austral. Math. Soc. Ser. A \textbf{31}, 202--206, (1981).
\end{thebibliography}


\end{document}